\documentclass[runningheads]{llncs}
\usepackage{graphicx}
\usepackage{amsmath, amssymb}
\usepackage{color}

\begin{document}
\title{Reduced Sum Implementation of the BURA Method for Spectral Fractional
  Diffusion Problems}
\titlerunning{Reduced Sum Implementation of the BURA Method}
%
\author{Stanislav Harizanov \and
 Nikola Kosturski \and Ivan Lirkov \and Svetozar Margenov \and Yavor Vutov}
\authorrunning{S. Harizanov et al.}
\institute{Institute of Information and Communication Technologies, \\
Bulgarian Academy of Sciences, Sofia, Bulgaria \\
\email{sharizanov@parallel.bas.bg}
\email{kosturski@parallel.bas.bg}
\email{ivan.lirkov@iict.bas.bg}
\email{margenov@parallel.bas.bg}
\email{yavor@parallel.bas.bg}}
\maketitle              
\begin{abstract}

The numerical solution of spectral fractional diffusion problems in the form
${\mathcal A}^\alpha u = f$ is studied, where $\mathcal A$ is a selfadjoint elliptic operator in a bounded domain $\Omega\subset {\mathbb R}^d$, and $\alpha \in (0,1]$.
The finite difference approximation of the problem leads to the system ${\mathbb A}^\alpha {\mathbf u} =    {\mathbf f}$, where ${\mathbb A}$ is a sparse, symmetric and positive definite (SPD) matrix, and ${\mathbb A}^\alpha$ is defined by its spectral decomposition. In the case of finite element approximation, ${\mathbb A}$ is SPD with respect to the dot product associated with the mass matrix. The BURA method is introduced by the best uniform rational approximation of degree $k$ of $t^{\alpha}$ in $[0,1]$, denoted by $r_{\alpha,k}$.
Then the approximation ${\bf u}_k\approx {\bf u}$ has the form
${\bf u}_k = c_0 {\mathbf f} +\sum_{i=1}^k c_i({\mathbb A} - {\widetilde{d}}_i
{\mathbb I})^{-1}{\mathbf f}$, ${\widetilde{d}}_i<0$, thus requiring the solving of
$k$ auxiliary linear systems with sparse SPD
matrices. The BURA method has almost optimal computational complexity,
assuming that an optimal PCG iterative solution method is applied to the
involved auxiliary linear systems. The presented analysis shows that
the absolute values of first 
$\left\{{\widetilde{d}}_i\right\}_{i=1}^{k'}$ 
can be extremely large. 
In such a case the condition
number of ${\mathbb A} - {\widetilde{d}}_i {\mathbb I}$ is practically equal to
one. Obviously, such systems do not need preconditioning. The next question is
if we can replace their solution by directly multiplying ${\mathbf f}$ with
$-c_i/{\widetilde{d}}_i$. Comparative analysis of numerical results is presented
as a proof-of-concept for the proposed RS-BURA method.
\end{abstract}
\section{Introduction}\label{sec-1}
We consider the second order elliptic operator ${\mathcal A}$ in the bounded domain $\Omega\subset {\mathbb R}^d$, $d\in\{1,2,3\}$, assuming homogeneous boundary conditions on $\partial\Omega$. Let $\{\lambda_j,\psi_j\}_{j=1}^\infty$ be the eigenvalues and eigenfunctions of $\mathcal{A}$, and let $(\cdot,\cdot)$ stand for the $L^2$ dot product. The spectral fractional diffusion problem ${\mathcal A}^\alpha u = f$ is defined by the equality
\begin{equation}\label{eq1}
\mathcal{A}^\alpha u =
\sum_{j=1}^\infty \lambda_j^\alpha (u,\psi_j)\psi_j,
~~\mbox{and therefore} ~~
u = \sum_{j=1}^\infty \lambda_j^{-\alpha} (f,\psi_j)\psi_j.
\end{equation}

Now, let a $(2d+1)$-point stencil on a uniform mesh be used to get the finite difference (FDM) approximation of $\mathcal A$. Then, the FDM numerical solution of (\ref{eq1}) is given by the linear system
\begin{equation}\label{eq2}
\mathbb{A}^\alpha {\bf u} = {\bf f},
\end{equation}
where ${\bf u}$ and ${\bf f}$ are the related mesh-point vectors, $\mathbb{A}\in\mathbb{R}^{N\times N}$ is SPD matrix, and $\mathbb{A}^\alpha$ is defined similarly to (\ref{eq1}), using the spectrum $\{\lambda_{j,h},\Psi_{j,h}\}_{j=1}^N$ of $\mathbb{A}$.

Alternatively, the finite element method (FEM) can be applied for the
numerical solution of the fractional diffusion problem, if $\Omega$ is a
general domain and some unstructured (say, triangular or tetrahedral)
mesh is used. Then, $\mathbb{A}=\mathbb{M}^{-1}\mathbb{K}$, where $\mathbb{K}$
and $\mathbb{M}$ are the stiffness and mass matrices respectively, and
$\mathbb{A}$ is SPD with respect to the energy dot
product associated with $\mathbb{M}$.

Rigorous error estimates for the linear FEM approximation of (\ref{eq1}) are
presented in \cite{BP-15}. More recently, the mass lumping case is
analyzed in \cite{HLMMP-20}. The general result is that the relative accuracy
in $L^2$ behaves essentially as $O(h^{2\alpha})$ for both linear FEM or
$(2d+1)$-point stencil FDM discretizations.

A survey on the latest achievements in numerical methods for fractional diffusion problems is presented in \cite{HLM-20}. Although there are several different approaches
discussed there, all the derived algorithms can be interpreted as rational
approximations, see also \cite{H-20}. In this context, the advantages of the BURA
(best uniform rational approximation) method are reported. The BURA 
method is originally introduced in \cite{HLMMV-18}, see \cite{HLMMP-20} and
the references therein for some further developments. The method is
generalized to the case of Neumann boundary conditions in \cite{HKMV-20}. The
present study is focused on the efficient implementation of the method for
large-scale problems, where preconditioned conjugate gradient (PCG) iterative
solver of optimal complexity is applied to the arising auxiliary sparse SPD
linear systems.

The rest of the paper is organized as follows. The construction and some basic
properties of the BURA method are presented in Section \ref{sec-2}. The
numerical stability of computing the BURA parameters is discussed in Section
\ref{sec-3}. The next section is devoted to the question how large can the
BURA coefficients be, followed by the introduction of the new RS-BURA method based on reduced sum
implementation in Section \ref{sec-5}. A comparative analysis of
the accuracy is provided for the considered test problem ($\alpha=0.25$,
$k=85$ and reduced sum of $46$ terms), varying the parameter $\delta$,
corresponding to the spectral condition number of $\mathbb A$. Brief
concluding remarks are given at the end.

\section{The BURA Method}\label{sec-2}
Let us consider the min-max problem
\begin{equation}\label{eq3}
\min_{r_k(t)\in \mathcal{R}(k,k)} \max_{t\in[0,1]} |t^\alpha -
r_k(t)|=:E_{\alpha,k}, ~~\alpha\in(0,1),
\end{equation}
where $r_k(t) = P_k(t)/Q_k(t)$, $P_k$ and $Q_k$ are polynomials of degree $k$,
and $E_{\alpha,k}$ stands for the error of the $k$-BURA element $r_{\alpha,k}$. Following \cite{HLMMP-20} we introduce 
\begin{equation}\label{eq4}
\mathbb{A}^{-\alpha} \approx \lambda_{1,h}^{-\alpha}r_{\alpha,k}(\lambda_{1,h}\mathbb{A}^{-1}), ~~{\mbox {and respectively}} ~~{\bf u}_k = \lambda_{1,h}^{-\alpha}r_{\alpha,k}(\lambda_{1,h}\mathbb{A}^{-1}) {\bf f},
\end{equation}
where ${\bf u}_k$ is the BURA numerical solution of the fractional diffusion system (\ref{eq2}). The following relative error estimate holds true (see \cite{S-93,S-03} for more details)
\begin{equation}\label{eq5}
\frac{||{\bf u} - {\bf u}_k||_2}{||{\bf f}||_2} \le \lambda_{1,h}^{-\alpha} E_{\alpha,k} \approx \lambda_{1,h}^{-\alpha} 4^{\alpha+1} \sin \alpha\pi ~ e^{-2\pi\sqrt{\alpha k}} = O\left (e^{-\sqrt{\alpha k}}~ \right ).
\end{equation} 

Let us denote the roots of $P$ and $Q$ by $\xi_1, \ldots , \xi_k$ and $d_1, \ldots , d_k$ respectively. It is known that they interlace, satisfying the inequalities
\begin{equation}\label{eq6}
0>\xi_1>d_1>\xi_2>d_2>\ldots>\xi_k>d_k.
\end{equation}
Using (\ref{eq4}) and (\ref{eq6}) we write the BURA numerical solution ${\bf u}_k$ in the form
\begin{equation}\label{eq7}
{\bf u}_k = c_0 {\bf f} + \sum_{i=1}^k c_i ({\mathbb A} - \widetilde{d}_i {\mathbb I})^{-1} {\bf f},
\end{equation}
where ${\widetilde d}_i=1/d_i<0$ and $c_i>0$. Further $0>\widetilde{d}_k > \ldots > \widetilde{d}_1$ and the following property plays an important role in this study
\begin{equation}\label{eq8}
\lim_{k\rightarrow \infty}\widetilde{d}_1 = -\infty.
\end{equation}
The proof of (\ref{eq8}) is beyond the scope of the present article. Numerical data illustrating the behaviour of $\widetilde{d}_i$ are provided in Section \ref{sec-4}. 

\section{Computing the Best Uniform Rational Approximation}\label{sec-3}
Computing the BURA element $r_{\alpha,k}(t)$ is a challenging computational problem on its own. There are more than 30 years modern history of efforts in this direction. For example, in \cite{VC-92}
the computed data for $E_{\alpha,k}$ are reported for six values of
$\alpha\in (0,1)$ with degrees $k\le 30$ by using computer arithmetic with 200
significant digits. The modified Remez method was used for the derivation of $r_{\alpha,k}(t)$ in \cite{HLMMV-18}. The min-max problem
(\ref{eq3}) is highly nonlinear and the Remez algorithm is very sensitive to
the precision of the computer arithmetic. In particular, this is due to the
extreme clustering at zero of the points of alternance of the error function.  Thus the
practical application of the Remez algorithm is reduced to degrees $k$ up to
$10-15$ if up to quadruple-precision arithmetic is used.

The discussed difficulties are practically avoided by the recent achievements
in the development of highly efficient algorithms that exploit a representation of the rational approximant in barycentric form and apply greedy
selection of the support points. The advantages of this approach is
demonstrated by the method proposed in \cite{H-20}, which is based on AAA
(adaptive Antoulas-Anderson) approximation of $z^{-\alpha}$ for
$z\in[\lambda_{1,h},\lambda_{N,h}]$. Further progress in computational
stability has been reported in \cite{H-21} where the BRASIL algorithm is
presented. It is based on the assumption that the best rational approximation
must interpolate the function at a certain number of nodes $z_j$, iteratively
rescaling the intervals $(z_{j-1}, z_j)$ with the goal of equilibrating the
local errors.

In what follows, the  approximation $r_{\alpha,k}(t)$ used in BURA is computed utilizing the open access software implementation of BRASIL, developed by Hofreither \cite{BRASIL}. In Table~\ref{tab1} we show the minimal degree $k$ needed to get a certain targeted accuracy $E_{\alpha,k}\in (10^{-12},10^{-3})$ for $\alpha\in \{0.25,0.50,0.75,0.80,0.90 \}$.

\begin{table}[t]\centering
\renewcommand{\arraystretch}{1.2}
\caption{ \footnotesize Minimal degree $k$ needed to get the desired accuracy $E_{\alpha,k}$}
\label{tab1}
\begin{tabular}{|c|c|c|c|c|c|}
\hline
~~~$E_{\alpha,k}$~~~& $~~~\alpha$ = 0.25~~~&~~~$\alpha$ = 0.50~~~&~~~$\alpha$ = 0.75~~~&~~~$\alpha$ = 0.80~~~&~~~$\alpha$ = 0.90~~~\\ \hline
$10^{-3}$ & 7 & 4 & 3 & 3 & 2 \\
$10^{-4}$ & 12 & 7 & 4 & 4 & 3 \\
$10^{-5}$ & 17 & 9 & 6 & 6 & 5 \\
$10^{-6}$ & 24 & 13 & 9 & 8 & 7 \\
$10^{-7}$ & 31 & 17 & 11 & 11 & 9 \\
$10^{-8}$ & 40 & 21 & 14 & 13 & 11 \\
$10^{-9}$ & 50 & 26 & 18 & 16 & 14 \\
$10^{-10}$ & 61 & 32 & 21 & 20 & 17 \\
$10^{-11}$ & 72 & 38 & 26 & 24 & 20 \\
$10^{-12}$ & 85 & 45 & 30 & 28 & 24 \\
\hline
\end{tabular}
\end{table}
The presented data well illustrate the sharpness of the error estimate
(\ref{eq5}). The dependence on $\alpha$ is clearly visible.

It is impressive that for $\alpha=0.25$ we are able to get a relative
accuracy of order $O(10^{-12})$ for $k=85$. Now comes the question how
such accuracy can be realized (if possible at all) in real-life
applications. Let us note that the published test results for larger $ k $ are
only for one-dimensional test problems. They are considered representative due
to the fact that the error estimates of the BURA method do not depend on the
spatial dimension $ d\in\{1,2,3\}$.

However, as we will see, some new challenges and
opportunities appear if $d>1$ and iterative solution methods have to be used
to solve the auxiliary large-scale SPD linear systems that appear in
(\ref{eq7}). This is due to the fact that some of the coefficients in the
related sum are extremely large for large $k$.

\section{How Large Can the BURA Coefficients Be}\label{sec-4}
The BURA method has almost optimal computational complexity $O((\log N)^2N)$
\cite{HLMMP-20}. This holds true under the assumption that a proper PCG
iterative solver of optimal complexity is used for the systems with matrices
${\mathbb A} - \widetilde{d}_i {\mathbb I}$, $\widetilde{d}_i < 0$. For example, some
algebraic multigrid (AMG) preconditioner can be used as an efficient solver,
if the positive diagonal perturbation $-\widetilde{d}_i\mathbb{I}$ of ${\mathbb A}$ does not substantially alter its condition number. 

The extreme behaviour of the first part of the coefficients $\widetilde{d}_i$ is
illustrated in Table \ref{tab2}. The case $\alpha = 0.25$ and $k=85$ is
considered in the numerical tests. Based on (\ref{eq8}), the expectation is
that $-\widetilde{d}_1$ can be extremely large. In practice, the
displayed coefficients are very large for almost a half of
$-\widetilde{d}_i$, starting with $-\widetilde{d}_1 \approx 10^{45}$. The
coefficients $\widetilde{d}_i$ included in Table \ref{tab2} are 
representative, taking into account their monotonicity (\ref{eq6}). We observe decreasing monotonic behavior for the coefficients $c_i$ and the ratios $c_i/\widetilde{d}_i$. Recall that the former were always positive, while the latter -- always negative. The theoretical investigation of this observation is outside the scope of the paper.

\begin{table}[t]\centering
\renewcommand{\arraystretch}{1.2}
\caption{Behaviour of the coefficients $\widetilde{d}_i$ and $c_i$ in (\ref{eq7}) for $\alpha = 0.25$ and $k=85$}
\label{tab2}
\addtolength{\tabcolsep}{2pt}
\begin{tabular}{|c|c|c|} \hline
$\widetilde{d}_1      = -1.7789 \times 10^{45}$  & $c_1  = 1.4698 \times 10^{34}$ & 
$c_1/\widetilde{d}_1 = -8.2624\times 10^{-12}$   \\
$\widetilde{d}_2      = -5.0719 \times 10^{42}$  & $c_2  = 1.1593 \times 10^{32}$ & 
$c_2/\widetilde{d}_2 = -2.2858\times 10^{-11}$ \\
$\widetilde{d}_3      = -7.1427 \times 10^{40}$  & $c_3  = 3.7598 \times 10^{30}$ & 
$c_3/\widetilde{d}_3 = -5.2638\times 10^{-11}$ \\
$\widetilde{d}_4      = -2.1087 \times 10^{39}$  & $c_4  = 2.2875 \times 10^{29}$ & 
$c_4/\widetilde{d}_4 = -1.0848\times 10^{-10}$ \\
$\widetilde{d}_5      = -9.7799 \times 10^{37}$  & $c_5  = 2.0285 \times 10^{28}$ & 
$c_5/\widetilde{d}_5 = -2.0741\times 10^{-10}$  \\
$\vdots$ & $\vdots$ & $\vdots$ \\
$\widetilde{d}_{35}    = -4.7256 \times 10^{17}$ & $c_{35} = 4.3270 \times 10^{12}$ & 
$c_{35}/\widetilde{d}_{35} = -9.1565\times 10^{-6}$\\
$\widetilde{d}_{36}    = -1.6388 \times 10^{17}$ & $c_{36} = 1.9278 \times 10^{12}$ & 
$c_{36}/\widetilde{d}_{36} = -1.1764\times 10^{-5}$\\
$\widetilde{d}_{37}    = -5.7675 \times 10^{16}$ & $c_{37} = 8.6879 \times 10^{11}$ & 
$c_{37}/\widetilde{d}_{37} = -1.5064\times 10^{-5}$\\
$\widetilde{d}_{38}    = -2.0587 \times 10^{16}$ & $c_{38} = 3.9585 \times 10^{11}$ & 
$c_{38}/\widetilde{d}_{38} = -1.9228\times 10^{-5}$\\
$\widetilde{d}_{39}    = -7.4487 \times 10^{15}$ & $c_{39} = 1.8227 \times 10^{11}$ & 
$c_{39}/\widetilde{d}_{39} = -2.4470\times 10^{-5}$\\
\hline
\end{tabular}
\end{table}
The conclusion is, that for such large values of $-\widetilde{d}_i$ the condition
number of ${\mathbb A} - \widetilde{d}_i {\mathbb I}$ practically equals 1. Of
course, it is a complete nonsense to precondition such matrices. Solving
such systems may even not be the right approach. In the next section we will
propose an alternative.

\section{Reduced Sum Method RS-BURA}\label{sec-5}
Let $k'<k$ be chosen such that $|\widetilde{d}_i|$ is very large for $1\le i \le k-k'$. Then, the reduced sum method for implementation of BURA (RS-BURA) for larger $k$ is defined by the following approximation ${\bf u}\approx{\bf u}_{k,k'}$ of the solution of (\ref{eq2}):
\begin{equation}\label{eq9}
{\bf u}_{k,k'} = \left [c_0 - \sum_{i=1}^{k-k'} \frac{c_i}{\widetilde{d}_i} \right ] {\bf f} + \sum_{i=k-k'+1}^k c_i ({\mathbb A} - \widetilde{d}_i {\mathbb I})^{-1} {\bf f}.
\end{equation}
The RS-BURA method can be interpreted as a rational approximation of ${\mathbb A}^{-\alpha}$ of degree $k'$ denoted by ${\widetilde{r}}_{\alpha,k,k'}({\mathbb A})$, thus involving solves of only $k'$ auxiliary sparse SPD systems.

Having changed the variable $z:=t^{-1}$ we introduce the error indicator 
$$
{\widetilde {E}}_{\alpha,k,k'}^\delta=\max_{z\in[1,\delta^{-1}]} |z^{-\alpha} - {\widetilde{r}}_{\alpha,k,k'}(z)|.
$$

Similarly to (\ref{eq5}), one can get the relative error estimate
\begin{equation}\label{eq10}
\frac{||{\bf u} - {\bf u}_{k,k'}||}{||{\bf f}||}\le {\widetilde {E}}_{\alpha,k,k'}^\delta,
\end{equation}
under the assumption $1\le\lambda_{1,h}< \lambda_{N,h}\le 1/\delta$. Let us remember that for the considered second order elliptic (diffusion) operator $\mathcal A$, the extreme eigenvalues of $\mathbb A$ satisfy the asymptotic bounds $ \lambda_{1,h} = O(1) ~~ {\mbox {and} } ~~ \lambda_{N,h} = O(h^{-2})$.

Direct computations give rise to 
\begin{equation}\label{eq:error analysis}
{\widetilde{r}}_{\alpha,k,k'}(z)-r_{\alpha,k}(z)=-\sum_{i=1}^{k-k'}\left [\frac{c_i}{\widetilde{d}_i}+\frac{c_i}{z-\widetilde{d}_i}\right ]=\sum_{i=1}^{k-k'}\left(-\frac{c_i}{\widetilde{d}_i}\right)\frac{z}{z-\widetilde{d}_i}>0,
\end{equation}
and since $x/(x-y)$ is monotonically increasing with respect to both arguments, provided $x>0$ and $y<0$, we obtain
$$
{\widetilde {E}}_{\alpha,k,k'}^\delta-E_{\alpha,k}\le\frac{\delta^{-1}}{\delta^{-1}-\widetilde{d}_{k-k'}} ~\frac{(k-k')c_{k-k'}}{-\widetilde{d}_{k-k'}}.
$$
When $k'$ is chosen large enough and $\delta^{-1}\ll -\widetilde{d}_{k-k'}$, the orders of the coefficients $\{\widetilde{d}_i\}_1^{k-k'}$ differ (see Table~\ref{tab1}), which allows us to remove the factor $(k-k')$ from the above estimate. Therefore, in such cases the conducted numerical experiments suggest that the following sharper result holds true
\begin{equation}\label{eq: order}
ord\left({\widetilde {E}}_{\alpha,k,k'}^\delta-E_{\alpha,k}\right)\le ord(\delta^{-1})+ ord\left(c_{k-k'}/(-\widetilde{d}_{k-k'})\right)-ord(-\widetilde{d}_{k-k'}), 
\end{equation}
providing us with an a priori estimate on the minimal value of $k'$ that will guarantee that the order of ${\widetilde {E}}_{\alpha,k,k'}^\delta-E_{\alpha,k}$ is smaller than the order of $E_{\alpha,k}$. Again, the theoretical validation of \eqref{eq: order} is outside the scope of this paper.

\begin{figure}[htp]
\begin{tabular}{cc}                                                           
\includegraphics[width=0.49\textwidth]{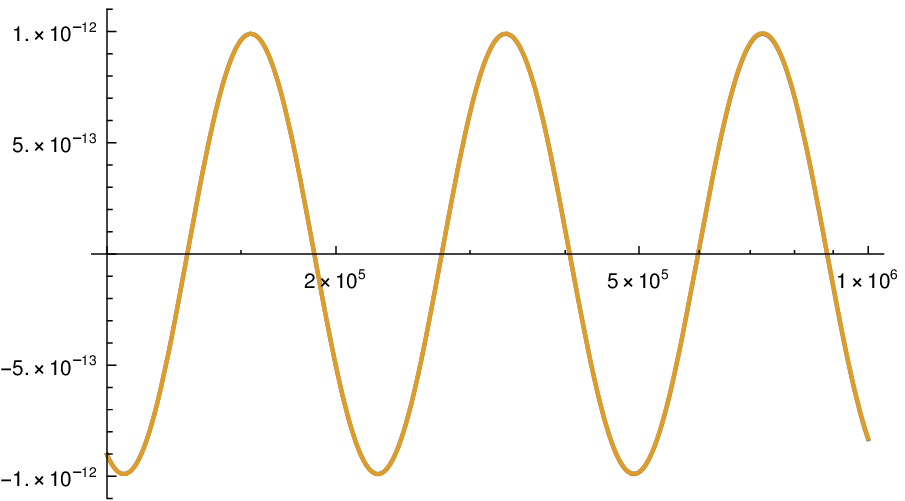}
\includegraphics[width=0.49\textwidth]{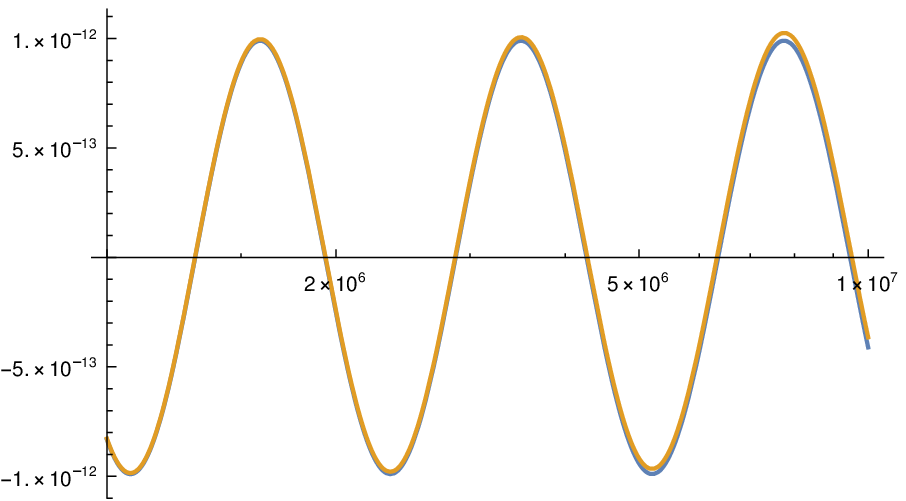}\\
\includegraphics[width=0.49\textwidth]{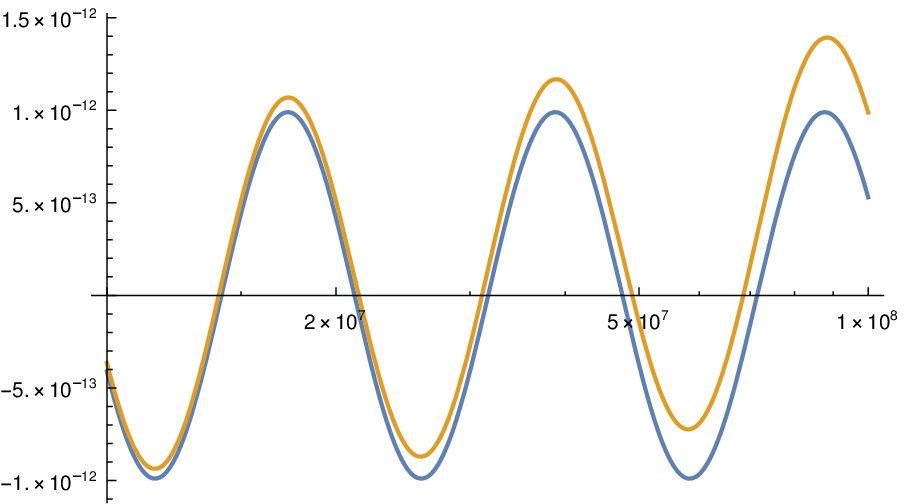}
\includegraphics[width=0.49\textwidth]{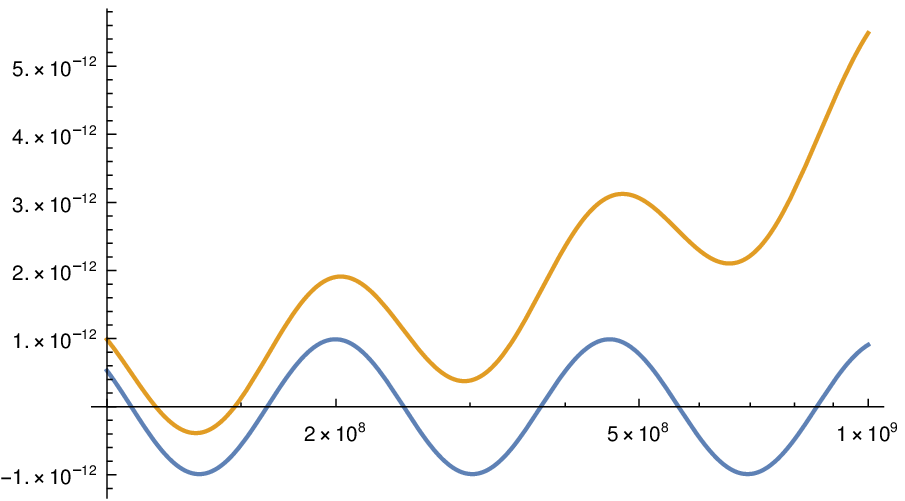}\\
\end{tabular}                                                                  
\vspace*{-0.5cm}
\caption{BURA (blue) and RS-BURA (yellow) errors $r_{\alpha,k}(z)-z^{-\alpha}$, respectively ${\widetilde{r}}_{\alpha,k,k'}(z)-z^{-\alpha}$, $z\in (1,\delta^{-1}]$, $\alpha=0.25$, $k=85$,
  $k'=46$, $\delta\in\{10^{-6},10^{-7},10^{-8},10^{-9}\}$.} 
\label{fig1}                               
\end{figure}

The next numerical results illustrate the accuracy of RS-BURA. The test
problem is defined by the following parameters: $\alpha=0.25$, $k=85$, $k'=46$
and $\delta\in\{10^{-6},10^{-7},10^{-8},10^{-9}\}$. In accordance with Table
\ref{tab1}, $E_{0.25,85}\approx 10^{-12}$. In Figure
\ref{fig1}, the behavior of  RS-BURA (yellow) and  BURA (blue) approximations
is compared. The RS-BURA error is always larger than the BURA error, as derived in \eqref{eq:error analysis}. The largest differences are located at the right part of the interval. As expected, for fixed $(k,k')$, the differences increase with the decrease of $\delta$. However, even in the worst case of $\delta=10^{-9}$, the accuracy ${\widetilde {E}}_{0.25,85,46}^\delta\approx 5
\times 10^{-12}$. Indeed, according to \eqref{eq: order} and Table~\ref{tab1}, in this case the order of ${\widetilde {E}}_{\alpha,k,k'}^\delta-E_{\alpha,k}$, thus the order of ${\widetilde {E}}_{\alpha,k,k'}^\delta$ itself, should not exceed $9+(-5)-15=-11$. On the other hand, \eqref{eq: order} implies that when $\delta>10^{-7}$ the order of ${\widetilde {E}}_{\alpha,k,k'}^\delta-E_{\alpha,k}$ is lower than the order of $E_{\alpha,k}$, meaning that the sum reduction does not practically affect the accuracy. This is illustrated on the first row of Figure~\ref{fig1}, where the plots of the two error functions are almost identical. This is a very promising result, taking into account that
$\delta$ is controlled by the condition number of $\mathbb A$.

From a practical point of view, if $\Omega \subset {\mathbb R}^3$, and
the FEM/FDM mesh is uniform or quasi uniform, the size of mesh parameter is
limited by $h > 10^{-4}$, corresponding to $N < 10^{12}$. Such
a restriction for $N$ holds true even for the nowadays
supercomputers. Thus, we obtain that if 
$h\approx 10^{-4}$, the corresponding  $\delta > 10^{-9}$. In general, for real
life applications, the FEM mesh is usually unstructured, including local
refinement. This leads to some additional increase of the condition number of
$\mathbb A$. However, the results in Figure \ref{fig1} show that RS-BURA has a
potential even for such applications.

The final conclusion of these tests is that some additional analysis is needed for better tuning of $k'$ for given $k$, $\alpha$ and a certain targeted accuracy.

\section{Concluding Remarks}\label{sec-6}
The BURA method was introduced in \cite{HLMMV-18}, where the number of the
auxiliary sparse SPD system solves was applied as a measure of computational
efficiency. This approach is commonly accepted in \cite{HLMMP-20,H-20},
following the unified view of the methods interpreted as rational
approximations, where the degree $k$ is used in the comparative analysis.

At the same time, it was noticed that when iterative PCG solvers are used in
BURA implementation even for relatively small degrees $k$, the number of
iterations is significantly different, depending on the values of
${\widetilde{d}}_i$ (see \cite{HKMV-20} and the references therein).

The present paper opens a new discussion about the implementation of BURA
methods in the case of larger degrees $k$. It was shown that due to the
extremely large values of part of  $|{\widetilde{d}}_i|$, the related auxiliary
sparse SPD systems are very difficult (if possible at all) to solve. In this
context, the first important question was whether BURA is even applicable to larger $k$ in
practice. The proposed RS-BURA method shows one possible
  way to do this. The presented numerical tests give promising indication for
the accuracy of the proposed reduced sum implementation of BURA.

The rigorous analysis of the new method (e.g., the proofs of \eqref{eq8}, \eqref{eq:error analysis}, and \eqref{eq: order}) is beyond the scope of this paper.

Also, the results presented in Section \ref{sec-5} reopen the topic of computational complexity analysis, including the question of whether the almost optimal estimate $O((\log N)^2N)$ can be improved. 

\subsection*{Acknowledgements}

We acknowledge the provided access to the e-infrastructure 
and support of the Centre for Advanced Computing and Data 
Processing, with the financial support by the Grant No 
BG05M2OP001-1.001-0003, financed by the Science and Education 
for Smart Growth Operational Program (2014-2020) and co-financed 
by the European Union through the European structural and 
Investment funds.

\vspace{2mm}
\noindent
The presented work is partially supported by the Bulgarian National Science Fund under grant No. DFNI-DN12/1.
\bibliography{LSSC21}
\end{document}